\documentstyle{article}
\begin{document}
\renewcommand{\baselinestretch}{1.2}
\begin{center}
{\Large\bf {Characterizations of all-derivable points in $B(H)$}{\footnote{This
work is supported by the Foundation of Hangzhou dianzi University
 }}} \vspace{0.4cm}
 \\
 { Jun Zhu{\footnote{E-mail address:
jzhu@hdu.edu.cn}}, Changping Xiong and Pan Li}\\
\vspace{0.3cm} \small{Institute of Mathematics, Hangzhou Dianzi
University, Hangzhou 310018, PR China\\}
\end{center}
\underline{~~~~~~~~~~~~~~~~~~~~~~~~~~~~~~~~~~~~~~~~~~~~~~~~~~~~~~~~~~~~~~~~~
~~~~~~~~~~~~~~~~~~~~~~~~~~~~~~~~~~~~~~~~~~~~~~~~~~~~~~~~~~~~~~~~~~~~~~}
\vspace{0.1cm}\\
{\bf{Abstract}} {

Let ${\mathcal{K}}$ and ${\mathcal{H}}$ be two Hilbert space, and let $B({\mathcal{K}},{\mathcal{H}})$ be the algebra of all
bounded linear operators from ${\mathcal{K}}$ into ${\mathcal{H}}$. We say that an element $G\in
B({\mathcal{H}},{\mathcal{H}})$ is an all-derivable point in $B({\mathcal{H}},{\mathcal{H}})$
if every derivable linear mapping $\varphi$ at $G$ (i.e.
$\varphi(ST)=\varphi(S)T+S\varphi(T)$ for any $S,T\in
B(H)$ with $ST=G$) is a derivation. Let both $\varphi: B({\mathcal{H}},{\mathcal{K}})\rightarrow B({\mathcal{H}},{\mathcal{K}})$ and $\psi: B({\mathcal{K}},{\mathcal{H}})\rightarrow B({\mathcal{K}},{\mathcal{H}})$ be two linear mappings. In this paper, the following results will be proved : if $Y\varphi(W)=\psi(Y)W$ for any $Y\in B({\mathcal{K}},{\mathcal{H}})$ and $W\in B({\mathcal{H}},{\mathcal{K}})$, then  $\varphi(W)=DW$ and $\psi(Y)=YD$ for some $D\in B({\mathcal{K}})$.
As an important application, we will show that an operator $G$ is an all-derivable point in
$B({\mathcal{H}},{\mathcal{H}})$ if and only if $G\neq 0$. \vspace{0.2cm}
\\
\vspace{0.2cm} {\it{AMS Classification}}: 47L30 15A30
47B47
\vspace{0.1cm}
\underline{~~~~~~~~~~~~~~~~~~~~~~~~~~~~~~~~~~~~~~~~~~~~~~~~~~~~~~~~~~~~~~~~~
~~~~~~~~~~~~~~~~~~~~~~~~~~~~~~~~~~~~~~~~~~~~~~~~~~~~~~~~~~~~~~~~~~~}

\section*{1. Introduction}
~

Let $\mathcal{A}$ be an algebra, and let
$\varphi$ be a linear mappings on ${\mathcal{A}}$. We say that
$\varphi$ is a derivation if $\varphi(ST)=\varphi(S)T+S\varphi(T)$
for any $S,T\in {\mathcal{A}}$. Fix an operator $G\in
{\mathcal{A}}$. We say that $\varphi$ is a derivable mapping at $G$
if $\varphi(ST)=\varphi(S)T+S\varphi(T)$ for any $S,T\in
 {\mathcal{A}}$ with $ST=G$. An element $G\in
 {\mathcal{A}}$ is called an all-derivable point in $ {\mathcal{A}}$
if every derivable mapping at $G$ is a derivation.

We describe some of the results related to ours.  Jin, Lu and Li
[2] show that every derivable mapping $\varphi$ at $0$ with
$\varphi(I)=0$ on nest algebras is a derivation. Hou and Qi [3]
prove that every derivable mapping at the unit operator on J-subspace
lattice algebras is a derivation. Zhu and Zhao [8] give the characterizations of all-derivable points in nest algebras $alg{\cal{N}}$ with nontrivial
nest $\cal{N}$(on Hilbert spaces), i.e. $G\in alg{\cal{N}}$ is an all-derivable point if and only if $G\neq 0$. But the condition whether the assumption of "nontrivial" may be omitted, remains open. It is obvious to see that a nontrivial nest algebra is essentially a triangular algebra, but trivial nest algebras $B({\mathcal{H}})$ is not triangular algebra, this case is more challenging than of triangular algebras.

Recently, Zhang, Hou and Qi [7] have proved the following theorem:

{\bf{Theorem 1.1 [7]}} {\it {Let $\cal{N}$ be a complete nest on a complex Banach space $\cal{X}$ with $dim {\cal{X}}\geq 2$, and $\delta: alg{\cal{N}}\rightarrow alg{\cal{N}}$ be a linear mapping. Let $Z\in alg{\cal{N}}$ be an injective operator or an operator with dense range in $alg{\cal{N}}$. Then $\delta$ is derivable at the operator $Z$ if and only if $\delta$ is a derivation. That is, every injective operator and every operator with dense range are all-derivable points of any nest algebra.}}

The above result implies the following corollary and it provides a basis for us to solve the above problem.

{\bf{Corollary 1.2 }} {\it{
Every injective operator and every operator with dense range are all-derivable points of a nest algebra( on Hilbert spaces).}}

The purpose of the present paper is to solve the above problem and prove that $G\in B({\mathcal{H}})$ is an all-derivable point if and only if $G\neq 0$. Furthermore, we obtain that $G\in alg{\cal{N}}$ is an all-derivable point if and only if $G\neq 0$ in any nest algebra.
For other results, see [1,4,5].

This paper is organized as follows. In Section 2, we describe the major Theorem 2.1 in
this paper, and give a preliminary Theorem 2.2 and its proof. Using the results, we give the proof of our main Theorem
2.1 in Section 3, that is $G$ is an all-derivable point if and only if $G\neq 0$ in any nest
algebras.

\section*{2. The main theorem and a lemma}
~

In this section we fix some notations. Let $\mathcal{H}$ and
$\mathcal{K}$ be two Hilbert spaces. We use the symbols $B({\mathcal{K}},{\mathcal{H}})$ and $F({\mathcal{K}},{\mathcal{H}})$
denote the set of all linear bounded operators from ${\mathcal{K}}$ into ${\mathcal{H}}$ and the set of all finite rank operators from ${\mathcal{K}}$ into ${\mathcal{H}}$, respectively.
The symbols $B({\mathcal{H}},{\mathcal{H}})$ and $F({\mathcal{H}},{\mathcal{H}})$ is
abbreviated to
$B({\mathcal{H}})$ and $F({\mathcal{H}})$, respectively.
Let $x\in {\mathcal{H}}$ and $y\in {\mathcal{K}}$, we use the symbols $x\otimes y\in B({\mathcal{K}},{\mathcal{H}})$ and $I_{\mathcal{H}}\in B({\mathcal{H}})$ denote the rank
one operator $<\cdot, y>x$ and the unit operator on ${\mathcal{H}}$,
respectively. Assume that $Y\in B({\mathcal{K}}, {\mathcal{H}})$, we always denote the range and kernel of $Y$ by the symbols $ran Y$ and $ker Y$, respectively.
we write $\mathbf{C}$ for the complex number field.

The following theorem is our main result:

{\bf{Theorem 2.1}} {\it Let ${\mathcal{H}}$ be a Hilbert space and
$G\in B(\mathcal{H})$.
Then $G$ is an all-derivable point in $B(\mathcal{H})$ if and only if $G\neq 0$.\\}

The following theorem will play an important role for our purposes.

{\bf{Theorem 2.2}} {\it Let ${\mathcal{H}}$ and ${\mathcal{K}}$ be two Hilbert spaces, and let $\varphi: B({\mathcal{H}}, {\mathcal{K}})\rightarrow B({\mathcal{H}}, {\mathcal{K}})$ and $\psi: B({\mathcal{K}}, {\mathcal{H}})\rightarrow B({\mathcal{K}}, {\mathcal{H}})$ be two linear mappings. If the following equation holds
\begin{eqnarray}
Y\varphi(W)=\psi(Y)W
\end{eqnarray}
for any $Y\in B({\mathcal{K}}, {\mathcal{H}})$ and $W\in B({\mathcal{H}}, {\mathcal{K}})$, then there exists an operator $D\in B({\mathcal{K}})$ such that $$\varphi(W)=DW~~~and~~~\psi(Y)=YD$$
for any $Y\in B({\mathcal{K}}, {\mathcal{H}})$ and $W\in B({\mathcal{H}}, {\mathcal{K}})$.}

{\bf Proof.} Step 1. Fix an operator $W\in B({\mathcal{H}}, {\mathcal{K}})$ with $ran W=\overline{ran W}^{\parallel\cdot \parallel}={\mathcal{K}}_{1}\subseteq{\mathcal{K}}$.
By the equation (1), we have
$$\psi(x\otimes y){{\mathcal{K}}_{1}}=\psi(x\otimes y)ran(W)\subseteq ran(x\otimes y)={\mathbf{C}}x$$
for any $x\in {\mathcal{H}}$ and $y\in {\mathcal{K}}$. Thus there exists a continuous linear functional $\lambda_{y,W}$ on ${\mathcal{K}}_{1}$ such that$$\psi(x\otimes y)(u)=\lambda_{y,W}(u)x$$
for any $u\in {\mathcal{K}}_{1}$. By Riesz representation theorem, there exists a vector $A_{W}(y)\in {\mathcal{K}}_{1}$ such that $\lambda_{y,W}(u)=<u, A_{W}(y)>$ and
$$\psi(x\otimes y)(u)=<u, A_{W}(y)>x=x\otimes A_{W}(y)(u)$$
for any $u\in {\mathcal{K}}_{1}$, i.e. $\psi(x\otimes y)=x\otimes A_{W}(y)$. By equation (1) and the above equation, we have
\begin{eqnarray}(x\otimes y)\varphi(W)=\psi(x\otimes y)W=x\otimes A_{W}(y)W.\end{eqnarray}
It is easy to verify from the linear of $\psi$ that $A_{W}:  {\mathcal{K}}\rightarrow  {\mathcal{K}}_{1}$ is a linear operator, and
$$\begin{array}{ccl}&&\mid<Wv, A_{W}(y)>\mid\parallel x\parallel=\parallel (x\otimes A_{W}(y))Wv\parallel=\parallel\psi(x\otimes y)Wv\parallel\\&&=\parallel x\otimes y\varphi(W)v\parallel\leq\parallel x\parallel\parallel y\parallel\parallel\varphi(W)v\parallel\end{array}$$
for any $y,v\in {\mathcal{K}}$. Fix a vector $u\in {\mathcal{K}}_{1}=ran W$, then the set $\{\mid<A_{W}(y),u>\mid: y\in {\mathcal{K}}_{1}, \parallel y\parallel\leq 1\}$ is a bounded set. By the uniform boundedness principle, $\{\parallel A_{W}(y)\parallel:  y\in {\mathcal{K}}_{1}, \parallel y\parallel\leq 1\}$  is a bounded set, i.e. $A_{W}\in B({\mathcal{K}},{\mathcal{K}}_{1})$.
It follows from the equation (2) that $(x\otimes y)\varphi(W)=x\otimes A_{W}(y)W=x\otimes y A_{W}^{*}W$, or
\begin{eqnarray}
\varphi(W)=A_{W}^{*}W.
\end{eqnarray}

Step 2. For any $x\in {\mathcal{H}}$ and $y\in {\mathcal{K}}$, then $y\otimes x\in B({\mathcal{H}}, {\mathcal{K}})$ and $ran (y\otimes x)=\overline{ran (y\otimes x)}^{\parallel\cdot\parallel}={\mathbf{C}}y$. It follows from equation (3) that $\varphi(y\otimes x)=A_{y\otimes x}(y\otimes x)$. Define a mapping $B_{y}: {\mathcal{H}}\rightarrow {\mathcal{H}}$ as follows:
$$B_{x}y=A_{y\otimes x}^{*}y$$
for any $y\in {\mathcal{K}}$. Then we have $$\varphi(y\otimes x)=B_{x}y\otimes x$$
for any $y\in {\mathcal{K}}$. It is easy to verify from the above equation that $B_{x}$ is a linear mapping on ${\mathcal{K}}$.
We claim that $B_{x}$ is independent of x. In fact, for any $x_{1},x_{2}\in {\mathcal{H}}$, we have
$$\varphi(y\otimes x_{i})=B_{x_{i}}y\otimes x_{i}~~(i=1,2);$$and
$$\varphi(y\otimes (x_{1}+x_{2}))=B_{x_{1}+x_{2}}y\otimes (x_{1}+x_{2}).$$
Combining the above three equations, we obtain
$$(B_{x_{1}}-B_{x_{1}+x_{2}})y\otimes x_{1}+  (B_{x_{2}}-B_{x_{1}+x_{2}})y\otimes x_{2}=0$$
If $x_{1}$ and $x_{2}$ are linearly independent, then $B_{x_{1}}y=B_{x_{1}+x_{2}}y=B_{x_{2}}y$ for any $y\in {\mathcal{K}}$, i.e. $B_{x_{1}}=B_{x_{2}}$. If $x_{2}=\alpha x_{1}$, then
$$B_{x_{2}}y\otimes x_{2}=\varphi(y\otimes x_{2})=\varphi(y\otimes \alpha x_{1})=\overline{\alpha}\varphi(y\otimes x_{1})=\overline{\alpha}B_{x_{1}}y\otimes x_{1}=B_{x_{1}}y\otimes \alpha x_{1}=B_{x_{1}}y\otimes x_{2}$$
for any $y\in {\mathcal{K}}$. Thus $B_{x_{2}}=B_{x_{1}}$. This implies that $B_{x}$ is independent of x. So we may write $D=B_{x}$ and get
$$\varphi(y\otimes x)=D(y\otimes x)$$
for any $y\in {\mathcal{K}}$ and $x\in {\mathcal{H}}$.
Hence \begin{eqnarray}\varphi(F)=DF\end{eqnarray} for any $F\in F({\mathcal{H}},{\mathcal{K}})$.

Step 3. Define a mapping $\Phi: B({\mathcal{H}},{\mathcal{K}})\rightarrow B({\mathcal{H}},{\mathcal{K}})$ as follows:
$$\Phi(W)=\varphi(W)-DW$$
for any $W\in B({\mathcal{H}},{\mathcal{K}})$.
We claim that $\Phi\equiv 0$. Otherwise, there exists $W\in B({\mathcal{H}},{\mathcal{K}})$ such that $\Phi(W)\neq 0$.
Furthermore we can find a vector $x\in {\mathcal{H}}$ with $\parallel x\parallel=1$ such that $y=\Phi(W)x\neq 0$.
Thus we have
$$\begin{array}{ccl}&&(x\otimes y)\Phi(W(I-x\otimes x))(x\otimes x)\\&=&(x\otimes y)[\varphi(W(I-x\otimes x))-DW(I-x\otimes x)](x\otimes x)\\&=&\psi(x\otimes y)W(I-x\otimes x)(x\otimes x)=0.\end{array}$$
Note that $\Phi(F)=0$ for any $F\in F({\mathcal{H}},{\mathcal{K}})$, so we have
$$\begin{array}{ccl}0&=&(x\otimes y)\Phi(W(I-x\otimes x))(x \otimes x)=(x\otimes y)[\Phi(W)-\Phi(Wx\otimes x)](x\otimes x)\\
&=&(x\otimes y)\Phi(W)(x\otimes x)=<\Phi(W)x,y>x\otimes x=\parallel y\parallel^{2}x\otimes x\end{array}$$
This is a contradiction with $x\neq 0$ and $y\neq 0$. Hence $\Phi\equiv0$, i.e. $\varphi(W)=DW$ for any $W\in B({\mathcal{H}},{\mathcal{K}})$.
Furthermore, by the equation (1), we obtain $\psi(Y)W =Y\varphi(W)=YDW$ for any $Y\in B({\mathcal{K}}, {\mathcal{H}})$ and $W\in B({\mathcal{H}},{\mathcal{K}})$.
Hence $$\psi(Y)=YD$$ for any $Y\in B({\mathcal{K}}, {\mathcal{H}})$.
 This completes the proof of the theorem.

 {\bf{Lemma 2.3.}} {\it Let $\mathcal{A}$ be an operator subalgebra
with unit operator $I$ in $B(H)$, and let $\varphi$ is a linear
mapping from $\mathcal{A}$ into $B({\mathcal{K}}_{1}, {\mathcal{K}}_{2})$. If $\varphi(X)=0$
for any invertible operator $X\in \mathcal{A}$, then $\varphi\equiv 0$.}

{\bf Proof.} For arbitrary operator $X\in {\mathcal{A}}$, there
exists a real number $\lambda>\parallel X\parallel$. Then both
$\lambda I-X$ and $2\lambda I-X$ are two invertible operators. So
$\varphi(\lambda I-X)=0$ and $\varphi(2\lambda I-X)=0$. It follows
from the linearity of $\varphi$ that $\varphi(X)=0$. Hence $\varphi\equiv 0$. This completes the proof of the Lemma.$\Box$

\section*{3. The proof of Theorem 2.1}
~

Now we will prove our main Theorem 2.1.

{\bf{The proof of Theorem 2.1.}} Necessary. Suppose that $G$ is an all-derivable point in $B({\mathcal{H}})$. We claim that $G\neq 0$.
In fact, the identity mapping $\varphi$ on $B({\mathcal{H}})$ is a derivable mapping at $0$, but $\varphi$ is not derivation.

Sufficiency.
Fix an operator $0\neq G\in B({\mathcal{H}})$ and $\varphi:
B(\mathcal{H})\rightarrow B(H)$ be a derivable mapping at
$G$.
We only need to prove that $\varphi$ is a derivation. If $\overline{ran G}={\mathcal{H}}$, then $\varphi$ is a derivation by Corollary 1.2.
If $\overline{ran W}\neq{\mathcal{H}}$, then we take ${\mathcal{H}}_{1}=\overline{ran G}$ and ${\mathcal{H}}_{2}=(ran G)^{\perp}$.
Obviously
$dim {\mathcal{H}}_{i}\geq 1 (i=1,2)$.
Then $G$ can be represented as a $2\times2$ operator matrix relative to the
orthogonal decomposition ${\mathcal{H}}={\mathcal{H}}_{1}\oplus {\mathcal{H}}_{2}$ as follows:
$$G=\left [\begin{array}{ccc}
  E & F \\
  0 & 0 \\
\end{array}\right ],$$
where $E\in B({\mathcal{H}}_{1})$, $F\in B({\mathcal{H}}_{2},{\mathcal{H}}_{1})$ with $E\neq 0$ or $F\neq 0$.
In the rest part of this paper, all $2\times2$ operator matrixes are always
represented as relative to the orthogonal decomposition ${\mathcal{H}}={\mathcal{H}}_{1}\oplus {\mathcal{H}}_{2}$.
For arbitrary operator $S\in  B(\mathcal{H})$, $S$ can be expressed as the following operator matrix
in the orthogonal decomposition of ${\mathcal{H}}={\mathcal{H}}_{1}\oplus {\mathcal{H}}_{2}$ as follows
$$S=\left [\begin{array}{ccc}
  X & Y \\
  Z & Q \\
\end{array}\right ],$$
where $X\in B({\mathcal{H}}_{1})$, $Y\in B({\mathcal{H}}_{2},{\mathcal{H}}_{1})$, $Z\in B({\mathcal{H}}_{1},{\mathcal{H}}_{2})$ and $Q\in B({\mathcal{H}}_{2})$. Since $\varphi$ is a linear mapping, we can write
$$\left\{\begin{array}{lll}\varphi(\left[\begin{array}{cc}
X & 0 \\
0 & 0
\end{array}\right])=\left[\begin{array}{cc}
A_{11}(X) & A_{12}(X)\\
A_{21}(X) &A_{22}(X)
\end{array}\right],\\~\\
\varphi(\left[\begin{array}{cc}
0 & Y \\
0 & 0
\end{array}\right])=\left[\begin{array}{cc}
B_{11}(Y) & B_{12}(Y)\\
B_{21}(Y) & B_{22}(Y)
\end{array}\right],\\~\\
\varphi(\left[\begin{array}{cc}
0 & 0 \\
Z & 0
\end{array}\right])=\left[\begin{array}{cc}
C_{11}(Z) & C_{12}(Z) \\
C_{21}(Z) & C_{22}(Z)
\end{array}\right],\\~\\
\varphi(\left[\begin{array}{cc}
0 & 0 \\
0 & Q
\end{array}\right])=\left[\begin{array}{cc}
D_{11}(Q) & D_{12}(Q) \\
D_{21}(Q) & D_{22}(Q)
\end{array}\right],\end{array}\right.$$
where $A_{ij}, B_{ij}, C_{ij}$ and $D_{ij}$ are linear mappings from ${\mathcal{H}}_{j}$ into ${\mathcal{H}}_{i}$.
For arbitrary $S,T\in B({\mathcal{H}})$ with $ST=G$, we may write
$$S=\left [\begin{array}{ccc}
  X & Y \\
  Z & Q \\
\end{array}\right ]$$
and
$$T=\left [\begin{array}{ccc}
  U & V \\
  W & R \\
\end{array}\right ].$$
It follows from $ST=G$ that $XU+YW=E, XV+YR=F, ZU+QW=0$ and $ZV+QR=0$. Since $\varphi$ is a derivable mapping at $G$,
we have
$$\begin{array}{ccl}&&\left[\begin{array}{ccc}
A_{11}(E)+B_{11}(F) & A_{12}(E)+B_{12}(F)\\
A_{21}(E)+A_{21}(F)] & A_{22}(E)+B_{22}(F)
\end{array}\right]=\varphi(G)=\varphi(S)T+S\varphi(T)\\&~&\\
&=&\left[\begin{array}{cc}
A_{11}(X)+B_{11}(Y)+C_{11}(Z)+D_{11}(Q) & A_{12}(X)+B_{12}(Y)+C_{12}(Z)+D_{12}(Q) \\
  A_{21}(X)+B_{21}(Y)+C_{21}(Z)+D_{21}(Q) &A_{22}(X)+B_{22}(Y)+C_{22}(Z)+D_{22}(Q)
\end{array}\right]\left[\begin{array}{cc}
U & V \\
W & R
\end{array}\right]\\&~&\\
&+&\left[\begin{array}{cc}
X & Y \\
Z & Q
\end{array}\right]\left[\begin{array}{cc}
A_{11}(U)+B_{11}(V)+C_{11}(W)+D_{11}(R) & A_{12}(U)+B_{12}(V)+C_{12}(W)+D_{12}(R)\\
A_{21}(U)+B_{21}(V)+C_{21}(W)+D_{21}(R) &A_{22}(U)+B_{22}(V)+C_{22}(W)+D_{22}(R)
\end{array}\right].
\end{array}$$
The above equation implies the following four operator equations
\begin{eqnarray}
&&A_{11}(E)+B_{11}(F)\nonumber\\
&=&A_{11}(X)U+B_{11}(Y)U+C_{11}(Z)U+D_{11}(Q)U\nonumber\\
&&+ A_{12}(X)W+B_{12}(Y)W+C_{12}(Z)W+D_{12}(Q)W\nonumber\\
&&+XA_{11}(U)+X B_{11}(V)+X C_{11}(W)+X D_{11}(R)\nonumber\\
&&+YA_{21}(U)+YB_{21}(V)+YC_{21}(W)+YD_{21}(R);
\end{eqnarray}
\begin{eqnarray}
&&A_{12}(E)+B_{12}(F)\nonumber\\
&=&A_{11}(X)V+B_{11}(Y)V+C_{11}(Z)V+D_{11}(Q)V\nonumber\\
&&+ A_{12}(X)R+B_{12}(Y)R+C_{12}(Z)R+D_{12}(Q)R\nonumber\\
&&+XA_{12}(U)+X B_{12}(V)+X C_{12}(W)+X D_{12}(R)\nonumber\\
&&+YA_{22}(U)+YB_{22}(V)+YC_{22}(W)+YD_{22}(R);
\end{eqnarray}
\begin{eqnarray}
&&A_{21}(E)+B_{21}(F)\nonumber\\
&=&A_{21}(X)U+B_{21}(Y)U+C_{21}(Z)U+D_{21}(Q)U\nonumber\\
&&+ A_{22}(X)W+B_{22}(Y)W+C_{22}(Z)W+D_{22}(Q)W\nonumber\\
&&+ZA_{11}(U)+ZB_{11}(V)+ZC_{11}(W)+ZD_{11}(R)\nonumber\\
&&+QA_{21}(U)+QB_{21}(V)+QC_{21}(W)+QD_{21}(R);
\end{eqnarray}
\begin{eqnarray}
&&A_{22}(E)+B_{22}(F)\nonumber\\
&=&A_{21}(X)V+B_{21}(Y)V+C_{21}(Z)V+D_{21}(Q)V\nonumber\\
&&+ A_{22}(X)R+B_{22}(Y)R+C_{22}(Z)R+D_{22}(Q)R\nonumber\\
&&+ZA_{12}(U)+ZB_{12}(V)+ZC_{12}(W)+ZD_{12}(R)\nonumber\\
&&+QA_{22}(U)+QB_{22}(V)+QC_{22}(W)+QD_{22}(R);
\end{eqnarray}
Note that the equations (5)-(8) always hold when $XU+YW=E, XV+YR=F, ZU+QW=0$ and $ZV+QR=0$.
Now we divide the proof of the theorem into the following eight Steps.

Step 1. We claim that $A_{22}\equiv 0$, $B_{21}\equiv 0$, $C_{12}\equiv 0$, $D_{22}\equiv 0$, $A_{11}(I_{{\mathcal{H}}_{1}})F=0$ and $A_{12}(YW)=YC_{22}(W)$ for any $W\in B({\mathcal{H}}_{1},{\mathcal{H}}_{2})$ and $Y\in B({\mathcal{H}}_{2},{\mathcal{H}}_{1})$.

1) For any invertible operator $X\in B({\mathcal{H}}_{1})$ and any $\lambda\in {\mathcal{C}}$, taking $Y=\lambda F$, $Z=0$, $Q=0$, $U=X^{-1}E$, $V=0$, $W=0$ and $R=\lambda^{-1}I_{{\mathcal{H}}_{2}}$ in the equation (8), then we have
$$A_{22}(E)+B_{22}(F)=\lambda^{-1}A_{22}(X)+B_{22}(F).$$
So $A_{22}(X)=0$. By Lemma 2.3, we obtain $A_{22}\equiv 0$.

2) For any operator $Y\in B({\mathcal{H}}_{2},{\mathcal{H}}_{1})$, taking $X=I_{{\mathcal{H}}_{1}}$, $Z=0$, $Q=0$, $U=E$, $V=F$, $W=0$ and $R=0$ in the equation (7) and (8) respectively, then we have
$$A_{21}(E)+B_{21}(F)=A_{21}(I_{{\mathcal{H}}_{2}})E+B_{21}(Y)E$$ and $$A_{22}(E)+B_{22}(F)=A_{21}(I_{{\mathcal{H}}_{2}})F+B_{21}(Y)F.$$
In the rest part of the proof, we always assume that $\lambda$ and $\mu$ are two arbitrary nonzero real number. Replacing $Y$ by $\lambda Y$ in the above two equations, we have
$$A_{21}(E)+B_{21}(F)=A_{21}(I_{{\mathcal{H}}_{2}})E+\lambda B_{21}(Y)E$$ and $$A_{22}(E)+B_{22}(F)=A_{21}(I_{{\mathcal{H}}_{2}})F+\lambda B_{21}(Y)F.$$
The above two equations implies that $B_{21}(Y)E=0$ and $B_{21}(Y)F=0$. It follows from $G=\left [\begin{array}{ccc}
  E & F \\
  0 & 0 \\
\end{array}\right ]$ that $B_{21}(Y)=B_{21}(Y)\mid_{{\mathcal{H}}_{2}}=B_{21}(Y)\mid_{ran G}=0$. Hence $B_{21}\equiv0$.

3) For any invertible operator $R\in B({\mathcal{H}}_{2})$, taking $X=I_{{\mathcal{H}}_{1}}$, $Y=FR^{-1}$, $Z=0$, $Q=0$, $U=E$, $V=0$ and $W=0$ in the equation (5), then we have
$$B_{11}(F)
=A_{11}(I_{{\mathcal{H}}_{1}})E+B_{11}(FR^{-1})E
+ D_{11}(R)
+FR^{-1}A_{21}(E)+FR^{-1}D_{21}(R).$$
Replacing $R$ by $\lambda R$ in the above equation, we have
$$B_{11}(F)
=A_{11}(I_{{\mathcal{H}}_{1}})E+FR^{-1}D_{21}(R)+\lambda^{-1}[B_{11}(FR^{-1})E+FR^{-1}A_{21}(E)]
+ \lambda D_{11}(R)
.$$
Hence we have
\begin{eqnarray}
B_{11}(F)=A_{11}(I_{{\mathcal{H}}_{1}})E+FR^{-1}D_{21}(R)
\end{eqnarray}
and $D_{11}(R)=0$. By Lemma 2.3, $D_{11}\equiv 0$.

4) For any $Y\in B({\mathcal{H}}_{2},{\mathcal{H}}_{1})$ and any $W\in B({\mathcal{H}}_{1},{\mathcal{H}}_{2})$, taking $X=I_{{\mathcal{H}}_{1}}$, $Z=0$, $Q=0$, $U=E-YW$, $V=F$ and $R=0$ in the equation (6), and noting that $A_{22}\equiv 0$, then we have
\begin{eqnarray}
0
=A_{11}(I_{{\mathcal{H}}_{1}})F+B_{11}(Y)F
-A_{12}(YW)+ C_{12}(W)
+YB_{22}(F)+YC_{22}(W).
\end{eqnarray}
Replacing $W$ by $\lambda W$ and taking $Y=0$ in the above equation (10), we have
$$0
=A_{11}(I_{{\mathcal{H}}_{1}})F+ \lambda C_{12}(W)
.$$
So $C_{12}(W)=0$ for any $W\in  B({\mathcal{H}}_{1},{\mathcal{H}}_{2})$, i.e. $C_{12}\equiv 0$.
On the other hand, replacing $W$ and $Y$ by $\lambda W$ and $\lambda Y$ in the equation (10), respectively, and Noting that $C_{12}\equiv 0$, then we have
$$0
=A_{11}(I_{{\mathcal{H}}_{1}})F+\lambda [B_{11}(Y)F+YB_{22}(F)]
+\lambda ^{2}[YC_{22}(W)-A_{12}(YW)]
.$$
Thus we have
\begin{eqnarray}
A_{11}(I_{{\mathcal{H}}_{1}})F=0
\end{eqnarray}and
\begin{eqnarray}
A_{12}(YW)=YC_{22}(W)
\end{eqnarray}
for any $W\in B({\mathcal{H}}_{1},{\mathcal{H}}_{2})$ and $Y\in  B({\mathcal{H}}_{2},{\mathcal{H}}_{1})$.

Step 2. We claim that $D_{22}$ is a derivation and $D_{12}(R)=-A_{12}(I_{{\mathcal{H}}_{1}})R$ for any $R\in B({\mathcal{H}}_{2})$.

For any $Y\in B({\mathcal{H}}_{2},{\mathcal{H}}_{1})$ and any $R\in B({\mathcal{H}}_{2})$, taking $X=I_{{\mathcal{H}}_{1}}$, $Z=0$, $Q=0$, $U=E$, $V=F-YR$ and $W=0$ in the equation (6), and using the results in Step 1, then we have
$$\begin{array}{ccl}
A_{12}(E)+B_{12}(F)
&=&-A_{11}(I_{{\mathcal{H}}_{1}})YR+B_{11}(Y)F-B_{11}(Y)YR\\&&
+ A_{12}(I_{{\mathcal{H}}_{1}})R+B_{12}(Y)R+A_{12}(E)
\\&&+ B_{12}(F)-B_{12}(YR)+ D_{12}(R)\\&&
+YB_{22}(F)-YB_{22}(YR)+YD_{22}(R).\end{array}
$$
Replacing $Y$ and $R$ by $\lambda Y$ and $\mu R$ in the above equation, respectively, then we obtain
$$\begin{array}{ccl}
0
&=&\lambda [B_{11}(Y)F+YB_{22}(F)]+\lambda \mu [YD_{22}(R)+ B_{12}(Y)R- A_{11}(I_{{\mathcal{H}}_{1}})YR\\&&- B_{12}(YR)]-\lambda ^{2}\mu [B_{11}(Y)YR- YB_{22}(YR)]+ \mu [A_{12}(I_{{\mathcal{H}}_{1}})R+ D_{12}(R)]
.\end{array}
$$
The above equation implies that
\begin{eqnarray}
B_{12}(YR)=B_{12}(Y)R+YD_{22}(R)-A_{11}(I_{{\mathcal{H}}_{1}})YR
\end{eqnarray}
and \begin{eqnarray}D_{12}(R)=-A_{12}(I_{{\mathcal{H}}_{1}})R\end{eqnarray}
for any $Y\in B({\mathcal{H}}_{2},{\mathcal{H}}_{1})$ and $R\in B({\mathcal{H}}_{2})$.

We claim that $A_{11}(I_{{\mathcal{H}}_{1}})=0$. In fact, for any $V\in Y\in B({\mathcal{H}}_{2},{\mathcal{H}}_{1})$,
taking $X=\lambda I_{{\mathcal{H}}_{1}}$, $Y=\mu (F-\lambda V)$, $Z=0$, $Q=0$, $U=\lambda ^{-1}E$, $V=F-YR$, $W=0$ and $R=\mu^{-1}I_{{\mathcal{H}}_{2}})$ in the equation (6), then we have
$$\begin{array}{ccl}0
&=&\lambda [A_{11}(I_{{\mathcal{H}}_{1}})V- VD_{22}(I_{{\mathcal{H}}_{2}})]+\mu [B_{11}(F)V+ FB_{22}(V)]\\
&&- \lambda \mu[ B_{11}(V)V+ VB_{22}(V)]
+ \lambda\mu^{-1}[A_{12}( I_{{\mathcal{H}}_{1}})
+ D_{12}(I_{{\mathcal{H}}_{2}}))]
\\&&+FD_{22}(I_{{\mathcal{H}}_{2}}).\end{array}$$
Thus we have
$$A_{11}(I_{{\mathcal{H}}_{1}})V= VD_{22}(I_{{\mathcal{H}}_{2}})$$
for any $V\in B({{\mathcal{H}}_{2}},{{\mathcal{H}}_{1}})$.
The above equation implies that $A_{11}(I_{{\mathcal{H}}_{1}})(x\otimes y)= (x\otimes y)D_{22}(I_{{\mathcal{H}}_{2}})$ for any $x\in {\mathcal{H}}_{1}$ and $y\in {\mathcal{H}}_{2}$. It follows that $A_{11}(I_{{\mathcal{H}}_{1}})x=\lambda_{x}x$. Hence $A_{11}(I_{{\mathcal{H}}_{1}})=\gamma I_{{\mathcal{H}}_{1}}$ for some $\gamma\in {\mathcal{C}}$.
We only need to prove that $\gamma=0$. Otherwise, $\gamma\neq 0$. By the equation (9) in Step 1, we have
$$B_{11}(F)=A_{11}(I_{{\mathcal{H}}_{1}})E+FR^{-1}D_{21}(R)=\gamma E+FR^{-1}D_{21}(R).$$
Supposing that $F=0$, it follows from the above equation that $E=0$. This is contradiction with $G\neq 0$. Supposing $F\neq 0$. Then $A_{11}(I_{{\mathcal{H}}_{1}})F=\gamma F\neq0$.
This is contradiction with the equation (11) in Step 1. Hence $A_{11}(I_{{\mathcal{H}}_{1}})=0$.

Now we show that $D_{22}$ is a derivation. For any $R_{1}, R_{2}\in B({{\mathcal{H}}_{2}})$ and $Y\in  B({{\mathcal{H}}_{2}},{{\mathcal{H}}_{1}})$, It follows from the equation (13) and $A_{11}(I_{{\mathcal{H}}_{1}})=0$ that
$$B_{12}(YR_{1}R_{2})=B_{12}(Y)R_{1}R_{2}+YD_{22}(R_{1}R_{2})$$
 and
$$B_{12}(YR_{1}R_{2})=B_{12}(Y R_{1})R_{2}+Y R_{1}D_{22}(R_{2})=B_{12}(Y) R_{1}R_{2}+YD_{22}(R_{1})R_{2}+Y R_{1}D_{22}(R_{2}).$$
The above two equations implies that $Y[D_{22}(R_{1}R_{2})-D_{22}(R_{1})R_{2}-R_{1}D_{22}(R_{2})]=0,$
i.e. $D_{22}( R_{1}R_{2})=D_{22}(R_{1})R_{2}-R_{1}D_{22}(R_{2})$ for any $R_{1},R_{2}\in B({{\mathcal{H}}_{2}})$.
Hence $D_{22}$ is a derivation. Since every derivation is an inner derivation on $B({\mathcal{H}}_{2})$, there exists an operator $\widetilde{D}\in B({\mathcal{H}}_{2})$ such that $D_{22}(R)=R\widetilde{D}-\widetilde{D}R$ for any $R\in B({{\mathcal{H}}_{2}})$.

Step 3. We claim that $A_{11}$ is a derivation and $B_{12}(XV)=X B_{12}(V)+ A_{11}(X)V$ for any operator $X\in B({\mathcal{H}}_{1})$ and $V\in B({\mathcal{H}}_{2},{\mathcal{H}}_{1})$.

For any invertible operator $X\in B({\mathcal{H}}_{1})$ and $V\in B({\mathcal{H}}_{2}, {\mathcal{H}}_{1})$, taking $Y=F-XV$, $Z=0$, $Q=0$, $U=X^{-1}E$, $W=0$ and $R=I_{{\mathcal{H}}_{2}}$ in the equation (6), then we have
$$\begin{array}{ccl}
A_{12}(E)+B_{12}(F)
&=&A_{11}(X)V+B_{11}(F)V-B_{11}(XV)V+ A_{12}(X)+B_{12}(F)\\
&&-B_{12}(XV)+XA_{12}(X^{-1}E)+X B_{12}(V)+X D_{12}(I_{{\mathcal{H}}_{2}})\\
&&+FB_{22}(V)-XVB_{22}(V)+(F-XV)D_{22}(I_{{\mathcal{H}}_{2}}).
\end{array}$$
Note that $D_{22}(I_{{\mathcal{H}}_{2}})=0$ by $D_{22}(R)=R\widetilde{D}-\widetilde{D}R$.  Replacing $X$ and $V$ by $\lambda X$ and $\mu V$ in the above equation, respectively, we have
$$\begin{array}{ccl}
A_{12}(E)
&=&\lambda \mu [A_{11}(X)V+ X B_{12}(V)-B_{12}(XV)]+\mu [B_{11}(F)V+ FB_{22}(V)]\\
&&-\lambda \mu^{2}[B_{11}(XV)V+XVB_{22}(V)]+ \lambda [A_{12}(X)+ X D_{12}(I_{{\mathcal{H}}_{2}})]+XA_{12}(X^{-1}E).
\end{array}$$
The above equation implies that
\begin{eqnarray}
B_{11}(XV)V+XVB_{22}(V)=0
\end{eqnarray} and
$$B_{12}(XV)=X B_{12}(V)+ A_{11}(X)V$$
for any invertible operator $X\in B({\mathcal{H}}_{1})$ and $V\in B({\mathcal{H}}_{2},{\mathcal{H}}_{1})$.
Fixing an operator $V\in B({\mathcal{H}}_{2},{\mathcal{H}}_{1})$, and define an linear mapping $\varphi: B({\mathcal{H}}_{1})\rightarrow B({\mathcal{H}}_{2},{\mathcal{H}}_{1})$ as follows
$$\varphi(X)=B_{12}(XV)-X B_{12}(V)- A_{11}(X)V$$
for any $X\in B({\mathcal{H}}_{1})$. Then $\varphi(X)=0$ for any invertible operator $X\in B({\mathcal{H}}_{1})$.
It follows from Lemma 2.3 that $\varphi\equiv 0$, i.e.
\begin{eqnarray}
B_{12}(XV)=X B_{12}(V)+ A_{11}(X)V
\end{eqnarray}
for any operator $X\in B({\mathcal{H}}_{1})$ and $V\in B({\mathcal{H}}_{2},{\mathcal{H}}_{1})$.

Now we show that $A_{11}$ is a derivation. For any $X_{1}, X_{2}\in B({\mathcal{H}}_{1})$ and $V\in B({\mathcal{H}}_{2},{\mathcal{H}}_{1})$, we have
$$B_{12}(X_{1}X_{2}V)=X_{1}X_{2}B_{12}(V)+A_{11}(X_{1}X_{2})V$$
and $$\begin{array}{ccl}&&B_{12}(X_{1}X_{2}V)=X_{1}B_{12}(X_{2}V)+A_{11}(X_{1})X_{2}V\\
&=&X_{1}X_{2}B_{12}(V)+X_{1}A_{11}(X_{2})V+A_{11}(X_{1})X_{2}V\end{array}$$
The above two equation implies that $[A_{11}(X_{1}X_{2})-X_{1}A_{11}(X_{2})-A_{11}(X_{1})X_{2}]V=0$. So $A_{11}$ is a derivation.
Since every derivation is an inner derivation on $B({\mathcal{H}}_{1})$, there exists an operator $A\in B({\mathcal{H}}_{1})$ such that $A_{11}(X)=XA-AX$ for any $X\in B({{\mathcal{H}}_{1}})$.

Step 4. We claim that $B_{22}(Y)=A_{21}(I_{{\mathcal{H}}_{1}})Y$ and $B_{11}(Y)=YD_{21}(I_{{\mathcal{H}}_{2}})$ for any $Y\in B({{\mathcal{H}}_{2}},{{\mathcal{H}}_{1}})$, and $D_{21}(R)=RD_{21}(I_{{\mathcal{H}}_{2}})$ for any $R\in B({{\mathcal{H}}_{2}})$.

1) For any $V\in B({{\mathcal{H}}_{2}},{{\mathcal{H}}_{1}})$, taking $X= I_{{\mathcal{H}}_{1}}$, $Y=F-V$, $Z=0$, $Q=0$, $U=E$, $W=0$ and $R=I_{{\mathcal{H}}_{2}}$ in the equation (8), then we have
\begin{eqnarray}
 B_{22}(V)=A_{21}(I_{{\mathcal{H}}_{1}})V
\end{eqnarray}
for any $V\in B({{\mathcal{H}}_{2}},{{\mathcal{H}}_{1}})$.

2) For any $Y\in B({{\mathcal{H}}_{2}},{{\mathcal{H}}_{1}})$ and $R\in B({{\mathcal{H}}_{2}})$, taking $X= \lambda I_{{\mathcal{H}}_{1}}$, $Z=0$, $Q=0$, $U=\lambda^{-1}E$, $V=\lambda^{-1}(F-YR)$ and $W=0$ in the equation (6), then we have
$$\begin{array}{ccl}
A_{12}(E)+B_{12}(F)
&=&B_{11}(Y)F-B_{11}(Y)YR
+\lambda A_{12}(I_{{\mathcal{H}}_{1}})R+B_{12}(Y)R\\
&&+ A_{12}(E)+ B_{12}(F)-B_{12}(YR)+\lambda D_{12}(R)
\\&&+\lambda^{-1} YB_{22}(F)-\lambda^{-1} YB_{22}(YR)+YD_{22}(R)
\end{array}$$
i.e.
$$\begin{array}{ccl}
0
=\lambda^{-1}[B_{11}(Y)E+YA_{21}(E)]
+[YD_{21}(R)-B_{11}(YR)]+\lambda D_{11}(R)
\end{array}
$$
The above equation implies that $$B_{11}(YR)=YD_{21}(R).$$ In particular, the following equation holds.
\begin{eqnarray}
 B_{11}(Y)=YD_{21}(I_{{\mathcal{H}}_{2}})
\end{eqnarray}
for any $Y\in B({{\mathcal{H}}_{2}},{{\mathcal{H}}_{1}})$. It follows that $YD_{21}(R)=B_{11}(YR)=YRD_{21}(I_{{\mathcal{H}}_{2}})$. Thus we have
\begin{eqnarray}
D_{21}(R)=RD_{21}(I_{{\mathcal{H}}_{2}})
\end{eqnarray}
for any $R\in B({{\mathcal{H}}_{2}})$.

Step 5. We claim that $C_{11}(W)=-A_{12}(I_{{\mathcal{H}}_{1}})W$ and $C_{22}(W)=WA_{12}(I_{{\mathcal{H}}_{1}})$ for any $W\in B({{\mathcal{H}}_{1}},{{\mathcal{H}}_{2}})$, and $A_{12}(X)=XA_{12}(I_{{\mathcal{H}}_{1}})$ for any $X\in B({{\mathcal{H}}_{1}})$.

For any $W\in B({{\mathcal{H}}_{1}},{{\mathcal{H}}_{2}})$ and $X\in B({{\mathcal{H}}_{1}})$, taking $Y=0$, $Z=0$, $Q=0$, $U=X^{-1}E$, $V=X^{-1}F$ and $R=0$ in the equation (5), then we have
$$
A_{11}(E)+B_{11}(F)
=A_{11}(X)X^{-1}E
+ A_{12}(X)W
+XA_{11}(X^{-1}E)+X B_{11}(X^{-1}F)+X C_{11}(W).
$$
Replacing $X$ by $\lambda X$ in the above equation, then we have
$$
A_{11}(E)+B_{11}(F)
=A_{11}(X)X^{-1}E
+XA_{11}(X^{-1}E)+X B_{11}(X^{-1}F)+ \lambda [A_{12}(X)W+ X C_{11}(W)].
$$
The above equation implies that $A_{12}(X)W=-X C_{11}(W)$. In particular,
 \begin{eqnarray}
C_{11}(W)=-A_{12}(I_{{\mathcal{H}}_{1}})W
\end{eqnarray}
for any $W\in B({{\mathcal{H}}_{1}},{{\mathcal{H}}_{2}})$. Furthermore
$
 A_{12}(X)W=-X C_{11}(W)=XA_{12}(I_{{\mathcal{H}}_{1}})W$. Hence
\begin{eqnarray}
A_{12}(X)=XA_{12}(I_{{\mathcal{H}}_{1}})
\end{eqnarray}
for any $X\in B({{\mathcal{H}}_{1}})$.
By the equation (12), we have
$$YC_{22}(W)=A_{12}(YW)=YWA_{12}(I_{{\mathcal{H}}_{1}})$$
The above equation implies that
\begin{eqnarray}
C_{22}(W)=WA_{12}(I_{{\mathcal{H}}_{1}})
\end{eqnarray}
for any $W\in B({{\mathcal{H}}_{1}},{{\mathcal{H}}_{2}})$.

Step 6. We claim that $A_{21}(X)=A_{21}(I_{{\mathcal{H}}_{1}})X$ for any $X\in B({{\mathcal{H}}_{1}})$.

For any invertible operator $X\in B({\mathcal{H}}_{1})$ and $V\in B({{\mathcal{H}}_{2}},{{\mathcal{H}}_{1}})$, taking $Y=F-XV$, $Z=0$, $Q=0$, $U=X^{-1}E$, $W=0$ and $R=I_{{\mathcal{H}}_{2}}$ in the equation (8), then we have
$$
B_{22}(XV)=A_{21}(X)V
$$
for any $X\in B({\mathcal{H}}_{1})$ and $V\in B({{\mathcal{H}}_{2}},{{\mathcal{H}}_{1}})$. It follows from the equation (17) that
$A_{21}(X)V=B_{22}(XV)=A_{21}(I_{{\mathcal{H}}_{1}})XV$. Thus we have
\begin{eqnarray}
A_{21}(X)=A_{21}(I_{{\mathcal{H}}_{1}})X
\end{eqnarray}
for any $X\in B({\mathcal{H}}_{1})$.

Step 7. We claim that $B_{12}(Y)=YD-AY$and $C_{21}(W)=WA-DW$ for some $D\in B({\mathcal{H}}_{2},{\mathcal{H}}_{1})$.

For any $Y\in  B({\mathcal{H}}_{2},{\mathcal{H}}_{1})$ and $W\in  B({\mathcal{H}}_{1},{\mathcal{H}}_{2})$,
taking $X=I_{{\mathcal{H}}_{1}}$, $Z=0$, $Q=0$, $U=E-YW$, $V=F$ and $R=0$ in the equation (5), then we have
$$\begin{array}{ccl}
0&=&
B_{11}(Y)E-B_{11}(Y)YW
+ A_{12}(I_{{\mathcal{H}}_{1}})W+B_{12}(Y)W-A_{11}(YW)
\\&&+ C_{11}(W)+YA_{21}(E)-YA_{21}(YW)+YB_{21}(F)+YC_{21}(W)
\end{array}$$
Replacing $Y$ and $W$ by $\lambda Y$ and $\lambda W$ in the above equation, we have
$$\begin{array}{ccl}
0&=&
\lambda [B_{11}(Y)E+ A_{12}(I_{{\mathcal{H}}_{1}})W+  C_{11}(W)+ YA_{21}(E)+ YB_{21}(F)]
\\&&+\lambda^{2}[B_{12}(Y)W-A_{11}(YW)+YC_{21}(W)]
-\lambda^{3}[YA_{21}(YW)+B_{11}(Y)YW].
\end{array}$$
The above equation implies that $$B_{12}(Y)W-A_{11}(YW)+YC_{21}(W)=0.$$ Since $A_{11}$ is a derivation and $A_{11}(X)=XA-AX$, the above equation implies that
$$
Y[WA-C_{21}(W)]=[B_{12}(Y)+AY]W
$$
for any $Y\in  B({\mathcal{H}}_{2},{\mathcal{H}}_{1})$ and $W\in  B({\mathcal{H}}_{1},{\mathcal{H}}_{2})$.
By Theorem 2.2, there exists $D\in B({\mathcal{H}}_{2})$ such that
\begin{eqnarray}
C_{21}(W)=WA-DW
\end{eqnarray} and
\begin{eqnarray}
B_{12}(Y)=YD-AY
\end{eqnarray}
for any $Y\in  B({\mathcal{H}}_{2},{\mathcal{H}}_{1})$ and $W\in  B({\mathcal{H}}_{1},{\mathcal{H}}_{2})$.
On the other hand,
for any $Y\in B({{\mathcal{H}}_{2}},{{\mathcal{H}}_{1}})$ and $R\in B({{\mathcal{H}}_{2}})$, taking $X= I_{{\mathcal{H}}_{1}}$, $Z=0$, $Q=0$, $U=E$, $V=F-YR$ and $W=0$ in the equation (6), then we have
$$\begin{array}{ccl}
0
&=&B_{11}(Y)F-B_{11}(Y)YR
+ A_{12}(I_{{\mathcal{H}}_{1}})R+B_{12}(Y)R\\
&&-B_{12}(YR)+ D_{12}(R)
+ YB_{22}(F)- YB_{22}(YR)+YD_{22}(R)
\end{array}$$
Replacing $Y$ and $R$ by $\lambda Y$ and $\lambda R$ in the above equation, we have
$$\begin{array}{ccl}
0
&=&\lambda[ B_{11}(Y)F+ A_{12}(I_{{\mathcal{H}}_{1}})R+  D_{12}(R)+  YB_{22}(F)]
\\
&&\lambda^{2}[B_{12}(Y)R- B_{12}(YR)
+ YD_{22}(R)]-\lambda^{3}[B_{11}(Y)YR+YB_{22}(YR)]
\end{array}$$
The above equation implies that $$B_{12}(YR)=B_{12}(Y)R
+ YD_{22}(R)$$
for any $Y\in B({{\mathcal{H}}_{2}},{{\mathcal{H}}_{1}})$ and $R\in B({{\mathcal{H}}_{2}})$.
It follows from the above equation and the equation (25) and $D_{22}(R)=R\widetilde{D}-\widetilde{D}R$ that
$$YRD-AYR=B_{12}(YR)=B_{12}(Y)R+YD_{22}(R)=YDR-AYR+YR\widetilde{D}-Y\widetilde{D}R,$$
i.e. $YR[D-\widetilde{D}]=Y[D-\widetilde{D}]R$. Furthermore $R[D-\widetilde{D}]=[D-\widetilde{D}]R$ for any $R\in B({{\mathcal{H}}_{2}})$.
Hence $D-\widetilde{D}=\alpha I_{{\mathcal{H}}_{2}}$ for some $\alpha\in {\mathcal{C}}$. It is easy to verify that
\begin{eqnarray}
D_{22}(R)=R\widetilde{D}-\widetilde{D}R=RD-DR
\end{eqnarray}
for any $R\in B({\mathcal{H}}_{2})$.

Step 8.
In summary, we have
$$\begin{array}{lll}
&A_{11}(X)=XA-AX,& B_{11}(y)=YD_{21}(I_{{\mathcal{H}}_{2}}), \nonumber\\
&A_{12}(X)=XA_{12}(I_{{\mathcal{H}}_{1}}),& B_{12}(y)=YD-AY,\nonumber\\
&A_{21}(X)=A_{21}(I_{{\mathcal{H}}_{1}})X,& B_{21}(y)=0,\nonumber\\
&A_{22}(X)=0,& B_{22}(y)=A_{21}(I_{{\mathcal{H}}_{1}})Y,\nonumber\\
&~&~\nonumber\\
&C_{11}(Z)=-A_{12}(I_{{\mathcal{H}}_{1}})Z,& D_{11}=0,\nonumber\\
&C_{12}(Z)=0,& D_{12}(Q)=-A_{12}(I_{{\mathcal{H}}_{1}})Q, \nonumber\\
&C_{21}(Z)=ZA-DZ,& D_{21}(Q)=QD_{21}(I_{{\mathcal{H}}_{2}}),\nonumber\\
&C_{22}(Z)=ZA_{12}(I_{{\mathcal{H}}_{1}}),&D_{22}=QD-DQ.
\end{array}$$
Combining the equation (15) with $B_{11}(Y)=YD_{21}(I_{{\mathcal{H}}_{2}})$ and $B_{22}(Y)=A_{21}(I_{{\mathcal{H}}_{1}})Y$, we have
$$0=B_{11}(XY)Y+XYB_{22}(Y)=XYD_{21}(I_{{\mathcal{H}}_{2}})Y+XYA_{21}(I_{{\mathcal{H}}_{1}})Y$$
for any $X\in B({{\mathcal{H}}_{1}})$ and $Y\in B({{\mathcal{H}}_{2}},{{\mathcal{H}}_{1}})$.
Hence we have $$D_{21}(I_{{\mathcal{H}}_{2}})=-A_{21}(I_{{\mathcal{H}}_{1}}).$$
For the convenience, we write $A_{12}(I_{{\mathcal{H}}_{1}})=B$ and $D_{21}(I_{{\mathcal{H}}_{2}})=C$.
Then we have
$$\left\{\begin{array}{lll}\varphi(\left[\begin{array}{cc}
X & 0 \\
0 & 0
\end{array}\right])=\left[\begin{array}{cc}
XA-AX & XB\\
-CX &0
\end{array}\right]=\left[\begin{array}{cc}
X & 0\\
0 &0
\end{array}\right]\left[\begin{array}{cc}
A & B\\
C & D
\end{array}\right]-\left[\begin{array}{cc}
A & B\\
C & D
\end{array}\right]\left[\begin{array}{cc}
X & 0\\
0 & 0
\end{array}\right],\\~\\
\varphi(\left[\begin{array}{cc}
0& Y \\
0 & 0
\end{array}\right])=\left[\begin{array}{cc}
YC & YD-AY\\
0 & -CY
\end{array}\right]=\left[\begin{array}{cc}
0 & Y\\
0 &0
\end{array}\right]\left[\begin{array}{cc}
A & B\\
C & D
\end{array}\right]-\left[\begin{array}{cc}
A & B\\
C & D
\end{array}\right]\left[\begin{array}{cc}
0 & Y\\
0 & 0
\end{array}\right],\\~\\
\varphi(\left[\begin{array}{cc}
0 & 0 \\
Z & 0
\end{array}\right])=\left[\begin{array}{cc}
-BZ & 0 \\
ZA -DZ & ZB
\end{array}\right]=\left[\begin{array}{cc}
0 & 0\\
Z &0
\end{array}\right]\left[\begin{array}{cc}
A & B\\
C & D
\end{array}\right]-\left[\begin{array}{cc}
A & B\\
C & D
\end{array}\right]\left[\begin{array}{cc}
0 & 0\\
Z & 0
\end{array}\right],\\~\\
\varphi(\left[\begin{array}{cc}
0 & 0 \\
0 & Q
\end{array}\right])=\left[\begin{array}{cc}
0 & -BQ \\
 QC  & QD-DQ
\end{array}\right]=\left[\begin{array}{cc}
0 & 0\\
0 & Q
\end{array}\right]\left[\begin{array}{cc}
A & B\\
C & D
\end{array}\right]-\left[\begin{array}{cc}
A & B\\
C & D
\end{array}\right]\left[\begin{array}{cc}
0 & 0\\
0 & Q
\end{array}\right]\end{array}\right.$$
for any $X\in B({\mathcal{H}}_{1})$, $Y\in B({\mathcal{H}}_{2},{\mathcal{H}}_{1})$, $Z\in B({\mathcal{H}}_{1},{\mathcal{H}}_{2})$ and $Q\in B({\mathcal{H}}_{2})$, i.e.
$$\varphi(\left[\begin{array}{cc}
X & Y\\
Z & Q
\end{array}\right])=\left[\begin{array}{cc}
X & Y\\
Z & Q
\end{array}\right]\left[\begin{array}{cc}
A & B\\
C & D
\end{array}\right]-\left[\begin{array}{cc}
A & B\\
C & D
\end{array}\right]\left[\begin{array}{cc}
X & Y\\
Z & Q
\end{array}\right]$$
for any $\left[\begin{array}{cc}
X & Y\\
Z & Q
\end{array}\right]\in B({\mathcal{H}})$. Hence $\varphi$ is an inner derivation.
This completes the proof. $\Box$

{\bf{Theorem 3.1}} {\it{Let ${\mathcal{N}}$ be a complete nest on Hilbert space ${\mathcal{H}}$. Then $G\in alg{{\mathcal{N}}}$ is an all-derivable point if and only if $G\neq 0$.}}

{\bf{ Proof.}} Supposing that ${\mathcal{N}}$ is a nontrivial complete nest, this is directly the conclusion of Theorem 2.3 in [7].
Supposing that ${\mathcal{N}}$ is a trivial complete nest, then $alg{\mathcal{N}}=B({\mathcal{H}})$. This is directly the conclusion of Theorem 2.1. This completes the proof. $\Box$
\section*{Reference}\small
\begin{description}
  \item[1] R.L. Crist, Local derivations on operator algebras, J. Funct. Anal. 135 (1996) 76-92.
  \item[2] W. Jing, S.J. Lu, P.T. Li,  Characterisations of
  derivations on some operator algebras, Bull. Austral. Math. Soc.
  66 (2002) 227-232.
  \item[3] J.C. Hou, X.F. Qi, Additive maps derivable at some points on J-subspace
  lattice algebras, Linear Algebra Appl. 429 (2008) 1851-1863.
  \item[4] D.R. Larson, A.R. Sourour, Local derivations and local automorphisms of ${\mathcal{B}}(X)$,
  Operator Algebras and Applications, Proc. Symp. Pure Math. 51
  (1990) 187-194.
  \item[5]  J.K. Li, Z.D. Pan, H. Xu,  Characterizations of isomorphisms and derivations of some algebras. J. Math. Anal. Appl.
  332 (2007) 1314-1322.
  \item[6] A.E. Taylor, D.C. Lay, Introduction to functional analysis(Second Edition), John Wiley and Sons, Inc. New York, 1980.
  \item[7] Y. F. Zhang, J.C. Hou, X.F. Qi, Characterizing derivations for any nest algebras on Banach space by their behaviors at an injective operator, Linear Algebra Appl., 449 (15) (2014) 312-333.
  \item[8]  J. Zhu, S. Zhao, Characterizations all-derivable points in nest algebras, Proc. Amer. Math. Soc. 141(7) (2013) 2343-2350.

  \end{description}

\end{document}